\title{ Constrained Ramsey Numbers}

\author{
Po-Shen Loh \thanks{Department of Mathematics,
Princeton University, Princeton, NJ 08544. E-mail:
 {\tt ploh@math.princeton.edu}.
Research supported in part by a Fannie and John Hertz Foundation Fellowship, an NSF Graduate
Research Fellowship, and a Princeton Centennial Fellowship.}
\and Benny Sudakov \thanks{
Department of Mathematics,
UCLA,  Los Angeles, CA 90095.
Email: {\tt bsudakov@math.ucla.edu}.
Research supported in part by NSF CAREER award DMS-0546523, NSF grant
DMS-0355497, USA-Israeli BSF grant, and by an Alfred P. Sloan fellowship.
}
}

\documentclass[11pt]{article}
\usepackage{amsmath, amssymb, amsthm, epsfig, graphicx, enumerate}

\newtheorem{theorem}{Theorem}[section]

\newtheorem{lemma}[theorem]{Lemma}
\newtheorem*{claim*}{Claim}

\newtheorem*{definition*}{Definition}

\newcommand{\E}[1]{\mathbb{E}\left[#1\right]}

\date{}
\begin{document}
\maketitle

\begin{abstract}
  For two graphs $S$ and $T$, the constrained Ramsey number $f(S, T)$ is the minimum 
  $n$ such that every edge coloring of the complete graph on $n$ vertices (with any number of
  colors) has a monochromatic subgraph isomorphic to $S$ or a rainbow subgraph isomorphic
  to $T$.  Here, a subgraph is said to be rainbow if all of its edges have different
  colors.  It is an immediate consequence of the Erd\H{o}s-Rado Canonical Ramsey Theorem
  that $f(S, T)$ exists if and only if $S$ is a star or $T$ is acyclic.  Much work has
  been done to determine the rate of growth of $f(S, T)$ for various types of parameters.
  When $S$ and $T$ are both trees having $s$ and $t$ edges respectively, Jamison, Jiang, and Ling
  showed that $f(S, T) \leq O(st^2)$ and conjectured that it is always at most $O(st)$.
  They also mentioned that one of the most interesting open special cases  is when $T$ is a 
  path. In this paper, we study this case and show that $f(S, P_t) = O(st\log t)$, 
  which differs only by a logarithmic factor from the conjecture. This substantially improves
  the previous bounds for most values of $s$ and $t$.
\end{abstract}

\section{Introduction}

The Erd\H{o}s-Rado Canonical Ramsey Theorem \cite{ER} guarantees that for any $m$, there
is some $n$ such that any edge coloring of the complete graph on the vertex set $\{1, \ldots,
n\}$, with arbitrarily many colors, has a complete subgraph of size $m$ whose coloring is
one of the following three types: monochromatic, rainbow, or lexical.
Here, a subgraph is rainbow if all edges receive distinct colors, and it 
is lexical when there is a total order of its vertices such that two 
edges have the same color if and only if they share the same larger endpoint.

Since the the first two types of colorings are somewhat more natural, 
it is interesting to study the cases when we can guarantee the existence of either 
monochromatic or rainbow subgraphs.
This motivates the notion of constrained Ramsey number $f(S, T)$, which is
defined to be the minimum $n$ such that every edge coloring of the complete graph on $n$
vertices (with any number of colors) has a monochromatic subgraph isomorphic to $S$ or a
rainbow subgraph isomorphic to $T$.  It is an immediate consequence of the Canonical
Ramsey Theorem that this number exists if and only if $S$ is a star or $T$ is acyclic,
because stars are the only graphs that admit a simultaneously lexical and monochromatic
coloring, and forests are the only graphs that admit a simultaneously lexical and rainbow
coloring.

The constrained Ramsey number has been studied by many researchers \cite{AJMP, BV, CSW,
 E-matchings, E-star-matching, GLS, JJL, JW, W}, and the bipartite case in \cite{BGLS}.
In the special case when $H=K_{1,k+1}$ is a star with $k+1$ edges, colorings with no 
rainbow $H$ have the property that every  vertex is incident to edges of at most $k$ different 
colors, and such colorings are called $k$-local. Hence  $f(S, K_{1,k+1})$
corresponds precisely to the local $k$-Ramsey numbers, $r^k_\text{loc} (S)$, which were
introduced and studied by Gy\'arf\'as, Lehel, Schelp, and Tuza in \cite{GLST}.
These numbers were shown to be within a
constant factor (depending only on $k$) of the classical $k$-colored Ramsey numbers $r(S;
k)$, by Truszczy\'nski and Tuza \cite{TT}.

When $S$ and $T$ are both trees having $s$ and $t$ edges respectively, Jamison, Jiang, and Ling
\cite{JJL} conjectured that $f(S, T) = O(st)$, and
provided a construction which showed that the conjecture, if true, is best
possible up to a multiplicative constant.
Here is a variant of such construction, which we present for the sake of completeness, which 
shows that in general the upper bound on $f(S, T)$ cannot be
brought below $(1+o(1))st$.  For a prime power $t$ let $\mathbb{F}_t$ be the finite field 
with $t$ elements. Consider the
complete graph with vertex set equal to the affine plane $\mathbb{F}_t
\times \mathbb{F}_t$, and color each edge based on
the slope of the line between the corresponding vertices in the affine
plane.  The number of different slopes (hence colors) is $t+1$, so
there is no rainbow graph with $t+2$ edges.  Also, monochromatic
connected components are cliques of order $t$, corresponding to affine lines. 
Therefore if $\Omega(\log t) < s < t$,  we can take a random subset of the
construction (taking each vertex independently with probability $s/t$) to obtain a coloring of the
complete graph of order $(1+o(1))st$ with $t+1$ colors in which 
all monochromatic connected components have size at most $(1+o(1))s$.  

Although Jamison, Jiang, and Ling were unable to prove their conjecture, 
they showed that $f(S, T) = O(st \cdot d_T) \leq O(st^2)$, where
$d_T$ is the diameter of $T$. Since this bound clearly gets weaker as the 
diameter of $T$ grows, they asked whether a pair of paths maximizes $f(S,T)$,
over all trees with $s$ and $t$ edges, respectively. This generated much 
interest in the special case when $T$ is a path $P_t$. In \cite{W}, 
Wagner proved that $f(S, P_t) \leq O(s^2 t)$. This bound grows linearly in $t$ when 
$s$ is fixed but still has order of magnitude $t^3$ for trees of similar size. Although
Gy\'arf\'as, Lehel, and Schelp \cite{GLS} recently showed that for small $t$ (less than 6),
paths are not the extremal example, they remain one of the most interesting
cases of the constrained Ramsey problem.

In this paper we prove the following theorem which agrees with the conjecture,
up to a logarithmic factor and the fact that $T$ is a path. It significantly improves
the previous bounds for most values of $s$ and $t$, and in particular gives the first 
sub-cubic bound for the case when the monochromatic tree and rainbow
path are of comparable size.

\begin{theorem}
  \label{thm:main}
  Let $S$ be any tree with $s$ edges, and let $t$ be a positive integer.  Then, for any $n
  \geq 3600st \log_2 t$, every coloring of the edges of the complete graph $K_n$ (with any
  number of colors) contains a monochromatic copy of $S$ or a rainbow
  $t$-edge path.
\end{theorem}

This supports the conjectured upper bound of $O(st)$ for the constrained Ramsey
number of a pair of trees.  With Oleg Pikhurko, the second author obtained another
result which provides further evidence for the conjecture. This result 
studies a natural relaxation of the above problem, in which one wants
to find  either a monochromatic copy of a tree $S$ or a properly colored copy of a tree $T$.
It appears that in this case the logarithmic factor can be removed, giving
an $O(st)$ upper bound.
We view this result as complementary to our main
theorem, and therefore have included its short proof in the appendix to our paper.

\vspace{4mm}

We close this section by comparing our approach to Wagner's, as the two proofs share some
similarities.  This will also lead us to introduce one of the the main tools that we will
use later.  Both proofs find a structured subgraph $G' \subset G$ in which one may direct
some edges in such a way that directed paths correspond to rainbow paths.  Wagner's
approach imposes more structure on $G'$, which simplifies the task of finding directed
paths, but this comes at the cost of substantially reducing $|G'|$.  In particular, his $|G'|$ is 
$s$ times smaller than $|G|$, which contributes a factor of $s$ to his ultimate bound
$O(s^2 t)$.  We instead construct a subgraph with weaker properties, but of order which is 
a constant fraction of $|G|$ (hence saving a factor of $s$ in the bound).  This
complicates the problem of finding the appropriate directed paths, which we overcome by
using the following notion of \emph{median order}:

\begin{definition*}
  Let $G$ be a graph, some of whose edges are directed.  Given a linear ordering $\sigma =
  (v_1, \ldots, v_n)$ of the vertex set, a directed edge $\overrightarrow{v_i v_j}$ is
  said to be \textbf{forward} if $i < j$, and \textbf{backward} if $i > j$.  If $\sigma$
  maximizes the number of forward edges, it is called a \textbf{median order}.
\end{definition*}

Median orders were originally studied for their own sake; for example, finding a median
order for a general digraph is known to be NP-hard. More recently,  Havet and Thomass\'e 
\cite{HT} discovered that they are a powerful tool for inductively building directed paths in
tournaments (complete graphs with all edges directed).  Their paper
used this method to produce a short proof of Dean's conjecture (see
\cite{DL}) that every tournament has a vertex whose second
neighborhood is at least as large as the first. Havet and Thomass\'e
also used a median order to attack Sumner's conjecture (see
\cite{W-Sumner}) that every tournament of order $2n-2$ contains every
oriented tree of order $n$. They succeeded in proving this conjecture
precisely for arborescences (oriented trees where every vertex except the
root has indegree one) and within a factor-2 approximation for
general oriented trees.

The only property that they used is the so-called \emph{feedback property}: if $\sigma =
(v_1, \ldots, v_n)$ is a median order, then for any pair $i < k$, the number of forward
edges $\overrightarrow{v_i v_j}$ with $i < j \leq k$ is at least the number of backward
edges $\overleftarrow{v_i v_j}$ with $i < j \leq k$.  This property is easily seen to be
true by comparing $\sigma$ to the linear order $\sigma' = (v_1, v_2, \ldots, v_{i-1},
v_{i+1}, v_{i+2}, \ldots, v_k, v_i, v_{k+1}, v_{k+2}, \ldots, v_n)$, which was obtained
from $\sigma$ by moving $v_i$ to the position between $v_k$ and $v_{k+1}$.
As an illustration of the simple power of this property, consider the following well-known
result, which we will in fact use later in our proof.

\begin{claim*}
   Every tournament has a directed Hamiltonian path.
\end{claim*}

\noindent {\bf Proof.}\, Let $\sigma = (v_1, \ldots, v_n)$ be a median
order.  For each $i$, the edge $v_i v_{i+1}$ is directed in
some way because we have a tournament, and so the feedback property
applied with $k=i+1$ implies that it is in fact a forward edge
$\overrightarrow{v_i v_{i+1}}$.  Therefore, $(v_1, \ldots, v_n)$ is
already a directed path, so we are done. \hfill $\Box$

\section{Proof of Theorem \ref{thm:main}}

Let us assume for the sake of contradiction that $n \geq 3600st \log_2 t$, but there is no
monochromatic copy of $S$ and no rainbow $t$-edge path.  In the past papers on 
the constrained Ramsey numbers of trees \cite{JJL, W}, and in this work, the following
well-known crude lemma is the only method used to exploit the nonexistence of a
monochromatic $S$. Its proof follows from the observation that every graph with average degree
$\geq 2s$ has an induced subgraph with minimum degree $\geq s$.

\begin{lemma}
  \label{lem:no-mono-edge-bound}
  Let $S$ be a tree with $s$ edges, and let $G = (V, E)$ be a simple graph, edge-colored
  with $k$ colors, with no monochromatic subgraph isomorphic to $S$.  Then $|E| < ks|V|$.
\end{lemma}

The rest of the proof of our main theorem roughly separates into two
main steps.  First, we find a structured subgraph $G' \subset G$ whose
order is within a constant factor of $|G|$.  We aim to
arrive at a contradiction by using $G'$ to construct a rainbow
$t$-edge path.  The structure of $G'$ allows us to direct many of its
edges in such a way that certain directed paths are automatically
rainbow.  In the second step, we use the median order's feedback
property to find many directed paths, which we then connect into a single long rainbow
path using the structure of $G'$.

\subsection{Passing to a directed graph}
\label{sec:structured-subgraph}

In this section, we show how to find a nicely structured subset of our
original graph, at a cost of a constant factor reduction of the size
of our vertex set.  We then show how the search for a rainbow path
reduces to a search for a particular collection of directed paths.

\begin{lemma}
  \label{lem:star-color-set} Let $S$ be a tree with $s$ edges and $t$
  be a positive integer.  Let $G$ be a complete graph on $n \geq
  310st$ vertices whose edges are colored (in any number of colors)
  in such a way that $G$ has no monochromatic copy of $S$ and no
  rainbow $t$-edge path.

  Then there exists a set $R$ of ``rogue colors'', a subset $U \subset V(G)$ with a
  partition $U = U_1 \cup \cdots \cup U_r$, an association of a distinct color $c_i \not \in
  R$ to each $U_i$, and an orientation of some of the edges of the induced subgraph $G[U]$, which satisfy the
  following properties:

  \begin{description}
  \item[(i)] $|U| > \frac{n}{10}$, $|R| < t$, and each $|U_i| < 2s$.
  \item[(ii)] For any edge between vertices $x \in U_i$ and $y \in U_j$ with $i \neq j$,
    if it is directed $\overrightarrow{xy}$, its color is $c_i$, if it is directed
    $\overrightarrow{yx}$, its color is $c_j$, and if it is undirected, its color is in $R$.
  \item[(iii)] For any pair of vertices $x \in U_i$ and $y \in U_j$ (where $i$ may equal $j$),
    there exist at least $t$ vertices $z \not \in U$ such that the
    color of the edge $xz$ is $c_i$ and the color of $yz$ is $c_j$.
  \end{description}
\end{lemma}

\noindent {\bf Proof.}\, Let us say that a vertex $v$ is
\emph{$t$-robust}\/ if for every set $F$ of $t$ colors, there are at
least $n/5$ edges adjacent to $v$ that are not in any of the colors in
$F$.  Let $V_1 \subset V$ be the set of $t$-robust vertices.  
We will need a lower bound on $|V_1|$, but this is just a special case
of Lemma \ref{lem:local-constraints} (whose short proof appears in the appendix).
Substituting the values $a=n/5$ and $b=t$ into this lemma gives 
$|V(G) \setminus V_1| \leq 2(ts+n/5) <4n/5$ and so $|V_1| \geq n/5$.

Now, let $P$ be a rainbow path of maximal length in $G$ such that at least one of its endpoints is
in $V_1$, and let $R$ be the set of colors of the edges of $P$.  $|R| < t$ by the
assumption that $G$ contains no rainbow $t$-edge path.  Let $B$ be the set of
vertices that have at least $n/15$ adjacent edges in a color in $R$.  Then $G$ contains at least $|B|n/30$ edges with colors in $R$.
On the other hand, by applying Lemma
\ref{lem:no-mono-edge-bound} to the subgraph of $G$ determined by taking only the edges
with colors in $R$, we see that the total number of edges in $G$ with color in $R$ is 
less than $|R| s n < tsn$, and so $|B| < 30st$.

Let $v$ be an endpoint of $P$ which is in $V_1$. Define the sets
$U_i$ as follows. Let $\{c_1, \ldots, c_r\}$ be the non-$R$ colors
that appear on edges adjacent to $v$. For each such $c_i$, let $U_i$ be the
set of vertices that are not in $B$ or $P$, and are adjacent to $v$ via an edge of color $c_i$.
Set $U = U_1 \cup \ldots \cup U_r$.  We
claim that these designations will satisfy the desired properties.

Consider arbitrary vertices $x \in U_i$ and $y \in U_j$, where $i$ may
equal $j$.  Since $n \geq 30t$, we have $|V_1 \setminus P| \geq
(2/15)n + t$, so $x, y \not \in B$ imply that there are at least $t$
choices for $z \in V_1 \setminus P$ such that both of the edges $xz$
and $yz$ have colors not in $R$.  Each such $xz$ must be in color
$c_i$, or else the extension of $P$ by the path $vxz$ would contradict
maximality of $P$, and similarly each $yz$ must be in color $c_j$.
Finally, $U \cap V_1 = \emptyset$, because any $w \in U \cap V_1$
would allow us to extend $P$ by the edge $vw$.  Therefore, we have
property (iii).

For property (ii), let $x \in U_i$ and $y \in U_j$, with $i \neq j$.  By property (iii),
there exists some vertex $z \in V_1 \setminus P$ such that $yz$ is in color $c_j$.  Then
the color of the edge $xy$ must be in $\{c_i, c_j\} \cup R$, or else the extension of $P$
by the path $vxyz$ would contradict its maximality.  Therefore, we can leave it undirected
if the color is in $R$, and direct it according to property (ii) otherwise.

It remains to show property (i).  We already established that $|R|<t$ and we can obtain the first inequality
from the construction of $V_1$ as follows.
Since $v \in V_1$, it is $t$-robust and so is adjacent to at least $n/5$ edges in non-$R$ colors.
Therefore, using that $n \geq 310st$ we get 
$$|U| \geq n/5-|B|-|P| > n/5- 30st - t \geq n/10.$$ 
For the last part, assume for
the sake of contradiction that $|U_i| \geq 2s$.  Arbitrarily select a
subset $U_i' \subset U_i$ of size $2s$, and consider the subgraph $G'$
formed by the edges of color $c_i$ among vertices in $U_i' \cup V_1$.  By the
argument that showed property (iii), every edge between $U_i'$ and
$V_1$ has color in $R \cup \{c_i\}$.  So, since $U_i \cap
B = \emptyset$, every $x \in U_i'$ is adjacent to at least $|V_1| -
n/15 \geq (2/3)|V_1|$ vertices in $V_1$ via edges of color $c_i$.
Therefore, using that $|V_1|/3 \geq 2s =|U_i'|$, we have
$$e(G') \geq |U_i'| \cdot (2/3)|V_1| = (4/3)s|V_1| =
s(|V_1| + (1/3)|V_1|) \geq s \cdot v(G').$$ 
Then Lemma \ref{lem:no-mono-edge-bound} implies that $G'$
has a copy of tree $S$, which is monochromatic by construction of
$G'$. This contradiction completes the proof of the last part of property (i), and the proof 
of the lemma. \hfill $\Box$

\vspace{4mm}

The partially directed subgraph of Lemma \ref{lem:star-color-set}
allows us to find rainbow paths by looking for certain types of
directed paths.  For example, if Lemma \ref{lem:star-color-set}
produces $U = U_1 \cup \ldots \cup U_m$, and we have found a directed
path $\overrightarrow{v_1 \ldots v_t}$ with each $v_i$ from a distinct
$U_j$, then it must be rainbow by property (ii) of the construction of $U$.  Unfortunately,
the following simple construction of a set with no monochromatic $S$
that satisfies the structure conditions of Lemma
\ref{lem:star-color-set} shows that we cannot hope to obtain our
rainbow path by searching for a single (long) directed path:
re-index $\{U_i\}$ with ordered pairs as
$\{U_{ij}\}_{i=1,j=1}^{h,t-1}$, let all $|U_{ij}| = s/3$, for all $1
\leq i < j \leq h$ direct all edges between any $U_{i,*}$ and
$U_{j, *}$ in the direction $U_{i, *} \rightarrow U_{j, *}$, 
and for all $1 \leq i \leq j < t$ and $1 \leq k \leq h$
color all edges between $U_{k, i}$ and $U_{k, j}$ in color $r_i$,
where $R = \{r_1, \ldots, r_{t-1}\}$.  Although it is clear that this
construction has no directed paths longer than $h =
O\big(\frac{|U|}{st}\big)$, it is also clear that one could build a
long rainbow path by combining undirected edges and directed paths.
The following lemma makes this precise.

\begin{lemma}
  \label{lem:linking}
  Let $U = U_1 \cup \ldots \cup U_m$ be a subset of $V(G)$ satisfying
  the structural conditions of Lemma \ref{lem:star-color-set}, and let
  $R$ be the associated set of rogue colors.  Suppose we have a
  collection of $r < t$ edges $\{u_iv_i\}_{i=1}^r$ in $G[U]$ whose colors are
  distinct members of $R$, and a collection of directed paths
  $\{P_i\}_{i=0}^r$, with $P_i$ starting at $v_i$ for $i \geq 1$.
  Then, as long as all of the vertices in $\{u_1, \ldots, u_r\} \cup
  P_0 \cup \ldots \cup P_r$ belong to distinct sets $U_j$, there exists a
  rainbow path in $G$ that contains all of the paths $P_i$ and all of
  the edges $u_iv_i$.  In short, one can link all of the fragments
  together into a single rainbow path.
\end{lemma}

\noindent {\bf Proof.}\, For each $i$, let $w_i$ be the final vertex
in the directed path $P_i$.  For a vertex $v \in U$, let $c(v)$ denote
the color associated with the set $U_i$ that contains $v$.  Since $r<t$, by
property (iii) of Lemma \ref{lem:star-color-set}, for each $0 \leq i < r$, there exists a distinct
vertex $x_i \not \in U$ such that the color of the edge $w_ix_i$ is
$c(w_i)$ and the color of the edge $x_iu_{i+1}$ is $c(u_{i+1})$. These vertices $x_i$ together with paths
$P_i$ form a
path $P$ of distinct vertices, which we will now prove is rainbow.

Note that our linking process only adds edges with non-rogue colors.
Since we assumed that the $u_iv_i$ have distinct colors, and the edges
of the $P_i$ are directed paths (hence with non-rogue colors), it is
immediate that $P$ has no duplicate rogue colors.  Also note that
among all directed edges in $\{P_i\}$, no pair of edges has initial
endpoint in the same $U_j$ by assumption.  Therefore, they all have
distinct colors by property (ii) of Lemma \ref{lem:star-color-set}.  Furthermore, none of these directed
edges originates from any point in any $U_j$ that intersects $\{u_1,
\ldots, u_r, w_1, \ldots, w_r\}$, so they share no colors with $C' =
\{c(u_1), \ldots, c(u_r), c(w_1), \ldots, c(w_r)\}$; finally the
colors in $C'$ are themselves distinct because of our assumption that
all vertices in $\{u_1, \ldots, u_r\} \cup P_0 \cup \ldots \cup P_r$
come from distinct $U_j$.  This proves that $P$ is a rainbow
path. \hfill $\Box$

\subsection{Finding directed paths}
\label{sec:build}

Now apply Lemma \ref{lem:star-color-set}, and let us focus on $U = U_1
\cup \ldots \cup U_m$, which is of size at least $n/10 \geq 360st
\log_2 t$.  Let us call the edges which have colors in $R$ ``rogue
edges.''  Note that if all edges were directed (i.e., we have a 
tournament), then the existence of a long directed
path follows from the fact that every tournament has a Hamiltonian
path.  The main issue is the presence of undirected edges.  We treat
these by observing that each undirected edge must have one of $|R|<t$ rogue
colors. Then, we use the machinery of median orders to
repeatedly halve the number of rogue colors, at the expense of losing
only $O(st)$ vertices each time.  This is roughly the source
of the $\log_2 t$ factor in our final bound.

Now we provide the details to make the above outline rigorous.  
Applying Lemma \ref{lem:no-mono-edge-bound} to the subgraph consisting
of all rogue edges, we see that the average rogue degree (number of
adjacent rogue edges) in $G[U]$ is at most $2s|R| \leq 2st$.  So, we
can delete all vertices in $U$ with rogue degree at least $4st$ at a
cost of reducing $|U|$ by at most half.  Let us also delete all edges
within each $U_i$ for the sake of clarity of presentation.  Note that
the reduced $U$ still has size at least $180st \log_2 t$.  Let
$\sigma$ be a median order for this partially directed graph induced
by $U$.  We will use the feedback property to find directed paths (and
this is the only property of median orders that we will use).

We wish to apply Lemma \ref{lem:linking}, so let us inductively build
a matching of distinct rogue colors, and accumulate a bad set that we
call $B$ and which we will maintain and update through the entire proof in this 
section.  Let $v_1$ be the first vertex according to $\sigma$, and
start with $B = U_\ell$, where $U_\ell \ni v_1$.  Proceed through the
rest of the vertices in the order of $\sigma$.  For the first stage,
stop when we first encounter a vertex not in $B$ that is adjacent to a
rogue edge (possibly several) whose other endpoint is also not in $B$, 
and call the vertex $v_2$.  Arbitrarily
select one of those rogue edges adjacent to $v_2$, call it $e_2$, and
call its color $r_2$.  Since we deleted all edges inside $U_i$, $e_2$
links two distinct $U_i$ and $U_j$.  Add all vertices of $U_i$ and
$U_j$ to $B$.  In general, if we already considered all vertices up to
$v_k$, continue along the median order (starting from the vertex
immediately after $v_k$) until we encounter a vertex not in $B$ that
is adjacent to an edge of a new rogue color which is not in $\{r_2,
\ldots, r_k\}$, again with other endpoint also not in $B$.  Call 
that vertex $v_{k+1}$, the edge $e_{k+1}$, and
its color $r_{k+1}$.  Add to $B$ all the vertices in the two sets
$U_i$ which contain the endpoints of $e_{k+1}$.  Repeat this procedure
until we have gone through all of the vertices in the order.  Suppose
that this process produces vertices $v_1, v_2,
\ldots, v_f$.  Then, to simplify the statements of our lemmas, also
let $v_{f+1}$, $v_{f+2}$, \ldots $v_{2f}$ refer to the final vertex in
the median order.  Our goal will be to find directed paths from
$\{v_i\}_{i=1}^f$, which via Lemma \ref{lem:linking}, will then extend
to a rainbow path.

Note that if $|B| \geq 2st$, then the number of vertices in $\{v_1\} \cup e_2 \cup \ldots
\cup e_f$ is at least $t$ by property (i).  Thus, applying Lemma \ref{lem:linking}  with
$P_i = \{v_i\}$, we can produce a rainbow path with at least $t$ edges.
Therefore, we may assume for the rest of this proof that $|B| < 2st$.  Also observe
that this argument implies that $f \leq t/2$.

The following technical lemma will help us to build the directed paths
$\{P_i\}$.

\begin{lemma}
  \label{lem:median-order-always-extend-if-st} Let $v$ be a vertex in
  $U$, and let $B$ be a set of size at most $2st$.  Then, among the
  $8st$ vertices immediately following $v$ in the median order, there is always
  some $w \not \in B$ such that there is a directed edge from $v$ to
  $w$.
\end{lemma}

\noindent {\bf Proof.}\, First, note that since we deleted all
vertices with rogue degree at least $4st$, more than $4st$ of the
$8st$ vertices immediately after $v$ are connected to $v$ by a
directed edge.  Since we have a median order, the feedback property
implies that only at most half of those edges can be directed back towards
$v$; therefore, there are more than $2st$ vertices there that have a
directed edge from $v$.  Since $|B| < 2st$, at least one of these
vertices will serve as our $w$.  \hfill $\Box$

\vspace{4mm}

Consider the vertices $v_1, v_2, v_4, \ldots, v_{2^{\lfloor \log_2 2f
\rfloor}}$.  Since we already established that $f \leq t/2$, this
is a list of at most $t+1$ vertices, the first and last
    of which are also the first and last vertices in the
    median order.  Since $U$ still has at least $180st\log_2 t$ vertices,
the pigeonhole principle implies that there must be
    some pair of vertices $\{v_\ell, v_{2\ell}\}$ in that list such
    that the number of vertices between them in the median order is at
    least $180st-2$.  Thus, the following lemma will provide the desired
    contradiction.

\begin{lemma}
  \label{lem:dyadic-build} If there is any $1 \leq \ell \leq f$ such
  that there are at least $176st$ vertices between $v_\ell$ and
  $v_{2\ell}$ in the median order, then $G$ has a rainbow $t$-edge
  path.
\end{lemma}

\noindent {\bf Proof.}\, Suppose we
have an $\ell$ that satisfies the conditions of the lemma.  Let $S_1$
be the first $8st$ vertices immediately following $v_\ell$ in the median
order, and let $S_2$ be the next $168st$ vertices in the median order.

Let us first build for every $i \leq \ell$ a
directed path $P_i$ from $v_i$  
to $S_1$ by repeatedly applying Lemma
\ref{lem:median-order-always-extend-if-st}. Start with each such $P_i
= \{v_i\}$, and as long as one of those $P_i$ does not reach $S_1$,
apply the lemma to extend it forward to a new vertex $w$, and add the
set $U_k$ containing $w$ to the set of bad vertices $B$.  If at any stage we have 
$|B| \geq 2st$,
we can immediately apply Lemma \ref{lem:linking} to find a rainbow
path with at least $t$ edges, just as in the argument directly
preceding the statement of Lemma
\ref{lem:median-order-always-extend-if-st}.  So, suppose that does not
happen, and let $\{w_i\}_1^\ell \subset S_1$ be the endpoints of these
paths.  We will show that we can further extend these paths into $S_2$
by a total amount of at least $t$, in such a way that we never use
two vertices from the same set $U_k$.  This will complete our proof
because Lemma \ref{lem:linking} can link them into a rainbow path
with at least $t$ edges.

Recall that all of the sets $U_i$ had size at most $2s$.  Therefore,
we can partition $S_2$ into  disjoint sets $U_j'$ with
$2s \leq |U_j'| \leq 4s$, where each $U_j'$ is obtained by taking a union of 
some sets $U_i \cap S_2$.  
We will design our path extension process such that it uses at most one vertex from
each $U_j'$, and hence it will also intersect each $U_k$ at most
once.  We use the probabilistic method to accomplish this.

Perform the following randomized algorithm, which will build a 
collection of sets $\{T_i\}_{i=1}^\ell$.  First, activate each $U_j'$ with
probability $1/8$.  Next, for each activated $U_j'$, select one of its
vertices uniformly at random, and assign it to one of the $T_i$, again
uniformly at random.  For each $i \leq \ell$, let $T_i'$ be
obtained from $T_i$ by deleting every vertex in $B$, and every vertex
that is not pointed to by a directed edge from $v_i$.  Finally, let $T_i''$ be
derived from $T_i'$ by (arbitrarily) deleting one vertex from every
rogue edge with both endpoints in $T_i'$.  Observe that now each
$T_i''$ spans a tournament, so as we saw at the end of the introduction, 
it contains a directed Hamiltonian path $P_i'$.  Since
$w_i$ has a directed edge to every vertex in $T_i''$, this $P_i'$ can
be used to extend $P_i$.  Therefore, if we can construct sets $T_i''$ such 
that $|T_1''| + \cdots +
|T_\ell''| \geq t$, we will be done.

Fix an $i \leq \ell$, and let us compute $\E{|T_i'|}$.  By the
feedback property of a median order, the number of (backward) directed
edges from $S_2$ to $\{w_i\}$ is at most half of the number of
directed edges between $w_i$ and the vertices in $S_1 \cup S_2$ which
follow it in the median order.  Since the latter number is bounded by
$|S_1 \cup S_2| = 176st$, the number of directed edges from $S_2$ to
$w_i$ is at most $88st$.  Also, the number of rogue edges between
$S_2$ and $\{w_i\}$ is at most $4st$ because by construction all rogue
degrees are bounded by $4st$.  Therefore, the number of (forward) directed edges
from $w_i$ to $S_2$ is at least $168st - 88st - 4st = 76st$.  Since we
will delete up to $2st$ vertices which were from $B$, the number of
directed edges from $w_i$ to vertices in $S_2 \setminus B$ is at least
$74st$.  Suppose $\overrightarrow{w_ix}$ is one of these directed
edges, and suppose that $x \in U_k'$.  The probability that $x$ is
selected for $T_i$ is precisely $\frac{1}{8 \cdot \ell \cdot |U_k'|}
\geq \frac{1}{8 \cdot \ell \cdot 4s}$, and by construction of $x$,
we know that if it is selected for $T_i$, it will also remain in
$T_i'$.  Therefore, by linearity of expectation,
\begin{displaymath}
  \E{|T_i'|} \geq 74st \cdot \frac{1}{8 \cdot \ell \cdot 4s} = \frac{37}{16} \frac{t}{\ell}
\end{displaymath}
To bound $\E{|T_i'| - |T_i''|}$, observe that the number of rogue
colors in the graph spanned by $S_2 \setminus B$ is less than $2\ell$,
by construction of the sequence $\{v_i\}$.  Therefore, Lemma
\ref{lem:no-mono-edge-bound} implies that there are less than $2\ell
\cdot s \cdot 168st$ rogue edges spanned by $S_2 \setminus B$.  Consider
one of these rogue edges $xy$.  
If we select both of its endpoints for $T_i$,
it will contribute at most 1 (possibly 0) to $|T_i'| - |T_i''|$;
otherwise it will contribute 0.  Above, we already explained
that the probability that the vertex $x \in U_j'$ is selected for
$T_i$ is precisely $\frac{1}{8 \cdot \ell \cdot |U_k'|} \leq
\frac{1}{8 \cdot \ell \cdot 2s}$.  If $x$ and $y$ come from distinct
$U_j'$, then the probabilities that they were both selected for $T_i$
are independent, and otherwise it is
impossible that they both were selected.  Hence
\begin{displaymath}
  \E{|T_i'| - |T_i''|} \leq 2\ell \cdot s \cdot 168st \cdot \left( \frac{1}{8 \cdot \ell \cdot 2s} \right)^2 = \frac{21}{16} \frac{t}{\ell}
\end{displaymath}
Therefore, by linearity of expectation, $\E{|T_i''|} \geq t/\ell$, and
thus $\E{|T_1''| + \cdots + |T_\ell''|} \geq t$.  This implies that 
there exists an instance of 
our random procedure for which $|T_1''| + \cdots + |T_\ell''| \geq t$, so we are 
done. \hfill 
$\Box$

\section{Concluding remarks}

In our proof, we apply Lemma \ref{lem:star-color-set} to produce a
structured set $U = U_1 \cup \ldots \cup U_m$ of size $\Omega(st \log
t)$.  The argument in Section \ref{sec:build} is quite wasteful
because, in particular, Lemma \ref{lem:dyadic-build} attempts to build
a collection of directed paths with total length $\geq t$, but
essentially using only the vertices in the median order between
$v_{\ell}$ and $v_{2 \ell}$.  This dissection of the vertex set into
dyadic chunks incurs the logarithmic factor in our bound.  We believe that
with a better argument, one might be able to complete the proof using
a structured set $U = U_1 \cup \ldots \cup U_m$ of size only
$\Omega(st)$.  If this were indeed possible, then Lemma
\ref{lem:star-color-set} would immediately imply that $f(S, P_t) =
O(st)$, because one loses only a constant factor in passing from
$V(G)$ to $U$.

It would be very interesting to obtain a better bound on $f(S, T)$ for
general trees $T$.  Our approach, based on the median order, seems
particularly promising here since it might be combined with the
following result of Havet and Thomass\'{e} \cite{HT} on Sumner's
conjecture: every tournament of order $4n$ contains every directed tree
of order $n$ as a subgraph.

\appendix

\section{Appendix (by Oleg Pikhurko and Benny Sudakov)}

Consider the following variant of the constrained Ramsey number. Let
$g(S,T)$ be the minimum integer $n$ such that every coloring of the
edges of the complete graph $K_n$ contains either a monochromatic copy
of $S$ or a properly colored copy of $T$. (In contrast, recall that
the definition of $f(S,T)$ requires $T$ to be rainbow). Similarly as for 
constrained Ramsey numbers, it is easy to see that $g(S, T)$ exists 
(i.e., it is finite) if and only if $S$ is a star or $T$ is acyclic. 
Although there has been little success bounding $f(S,T)$ by $O(st)$, it turns
out that we can prove a quadratic upper bound for $g(S,T)$, which is
of course no larger than $f(S,T)$.

\begin{theorem}
  \label{thm:properly-colored}
  Let $S$ and $T$ be two trees with $s$ and $t$ edges,
  respectively. Then $g(S,T) \leq 2st+t^2$.
\end{theorem}

The following construction shows that the upper bound is tight up to a
constant factor.  Let $S$ be a path with $s+1$ edges and $T$ be a star
with $t+1$ edges.  Then let $V_1, \ldots, V_t$ be disjoint sets of size
$\lfloor s/2 \rfloor$ each.  Color all edges inside $V_i$ and from $V_i$ to $V_j$ with
$j>i$ by color $i$.  This produces a graph on $t \lfloor s/2 \rfloor$ vertices with no
monochromatic $S$ and no properly colored $T$.

To prove Theorem \ref{thm:properly-colored}, we first need the
following lemma.

\begin{lemma}
\label{lem:local-constraints}
  Consider an edge coloring of the complete graph which contains no
  monochromatic copy of a fixed tree $S$ with $s$ edges. Let $U$ be
  the set of vertices such that for every $u \in U$, one can delete at
  most $a$ edges from the graph such that the remaining edges which
  connect $u$ to the rest of the graph have at most $b$ colors. Then
  $|U|\leq 2(bs+a)$.
\end{lemma}

\noindent {\bf Proof.}\, Focus on the subgraph induced by $U$. Now
we can delete at most $a$ edges at every vertex so that the
remaining edges at that vertex have at most $b$ colors. Let $G$ be the
graph obtained after all of these deletions. If $m = |U|$, then the
number of edges of $G$ is at least ${m \choose 2}-am$. For every
remaining color $c$, let $G_c$ be the subgraph of all edges of color
$c$.  By Lemma \ref{lem:no-mono-edge-bound}, we
have $e(G_c) < s \cdot v(G_c)$ for each $c$.  Also, since every vertex of $G$
is incident with edges of at most $b$ colors, we have that $\sum_c
v(G_c) \leq bm$.  Combining all these inequalities we have
\begin{displaymath}
  {m \choose 2}-am \leq \sum_c e(G_c) < \sum_c s \cdot v(G_c) \leq sbm.
\end{displaymath}
This implies that $m< 2(bs+a)+1$.
\hfill $\Box$

\vspace{4mm}

\noindent {\bf Proof of Theorem \ref{thm:properly-colored}.}\, The
proof is by induction on $t$.  The statement is trivial for $t=1$
because any edge will give us a properly colored $T$.  Now suppose
that $T$ is a tree with $t > 1$ edges, and we have a coloring of $G =
K_{2st+t^2}$ with no monochromatic copy of $S$.  It suffices to show
that we can find a properly colored copy of $T$.  Select an edge
$(u,v)$ of $T$ such that all neighbors of $v$ except $u$ are leaves
$v_1, \ldots, v_k$. Delete $v_1, \ldots, v_k$ from $T$ and call the
new tree $T_1$. The number of edges in $T_1$ is $t_1=t-k$.

Let $U$ be the set of vertices of $G$ such that for every $u \in U$
one can delete at most $t_1$ edges from $G$ such that the edges
which connect $u$ to the rest of the graph have at most $k$
colors. By the previous lemma $|U| \leq 2(ks+t_1)$, and let
$W=V(G) \setminus U$.  Then we have that
\begin{displaymath}
  |W| = 2st + t^2 - |U| \geq 2st_1 + t^2 - 2t_1 > 2st_1 + t_1^2.
\end{displaymath}
Therefore by induction we can find a properly colored copy of the tree
$T_1$ inside $W$.  Let $u',v'$ be the images in this copy of the
vertices of $u,v$ of $T_1$.  By definition of $W$, the vertex $v'$ has
edges of at least $k+1$ colors connecting it with vertices outside
this copy of $T_1$. At least $k$ of these colors are different from
that of the edge $(u',v')$, so we can extend the tree to a properly colored
copy of $T$. \hfill $\Box$

Using a more careful analysis in the above proof, which we omit, one can slightly improve 
the term $t^2$ in Theorem \ref{thm:properly-colored}.

\end{document}